\documentclass[a4paper,11pt]{amsart}
\usepackage[latin1]{inputenc}
\usepackage{amsmath,amsthm,amsfonts,amscd,amssymb,eucal,latexsym,mathrsfs}
\usepackage[all]{xy}
\usepackage{latexsym, amssymb, amscd, mathrsfs, graphics, fullpage,graphicx}
\usepackage[english,french]{babel}
\newtheorem{theorem}{Th\'eor\`eme}
\newtheorem{corollary}[theorem]{Corollaire}

\theoremstyle{definition}

\DeclareMathOperator{\dom}{dom}

\renewcommand{\Im}{\mathop{\mathrm{Im}}}

\newcommand{\ZI}{{\mathrm{\bf Z}}}

\newcommand{\R}{R}   
 
\newcommand{\G}{\Gamma}

\newcommand{\IIi}{\mathrm{II_1}}

\renewcommand{\d}{\mathrm{d}}

\newcommand{\del}{\partial}

\date{\today}

\title{Le co\^ut est un invariant isop\'erim\'etrique} 

\author{Mika\"el Pichot et St\'ephane Vassout}

\address{\hskip-\parindent
Mika\"el Pichot, IPMU,
University of Tokyo,
3-8-1 Komaba, Tokyo, 153-8914
JAPAN}
\email{pichot@ms.u-tokyo.ac.jp}

\address{\hskip-\parindent
St\'ephane Vassout, Institut de Math\'ematiques de Jussieu and Universit\'e Paris 7}
\email{vassout@math.jussieu.fr}

\begin{document}

\begin{abstract} Soit $R$ une relation d'\'equivalence mesur\'ee ergodique  de type $\IIi$ sur un espace de probabilit\'e $(X,\mu)$ sans atome. On montre que
\[
C(R)= 1+{1\over 2}h(R),
\]
o\`u   $C(R)$ est le co\^ut de $R$, suppos\'e fini, et  $h(R)$ sa constante isop\'erim\'etrique. Cela fait suite aux r\'esultats r\'ecents de Russell Lyons et des auteurs sur le sujet.
\end{abstract}

\newcommand{\val}{\mathrm{val}}

\maketitle

\section{Introduction}

Soit $R$ une  relation d'\'equivalence mesur\'ee ergodique  de type $\IIi$ (\`a classes d\'enombrables) sur un espace de probabilit\'e $(X,\mu)$ sans atome. Comme pour les groupes d\'enombrables, on associe \`a tout syst\`eme g\'en\'erateur $\Phi$ de $R$ un graphe de Cayley $\Sigma_\Phi$ sur lequel la relation $\R$ agit ($\Phi$ est aussi appel\'e un graphage de $R$). Suivant G. Levitt (voir \cite{Gaboriau99}), le \emph{co\^ut} de $R$ est d\'efini par 
\[
C(R)=\inf_{\Phi}\, {1\over 2} \val_\mu(\Sigma_\Phi)
\]
o\`u l'infimum porte sur les syst\`emes g\'en\'erateurs de $R$ et $\val_\mu$ est la valence moyenne dans $\Sigma_\Phi$ prise relativement \`a $\mu$. Il s'agit a priori d'un  nombre r\'eel compris entre 1 et $\infty$ (inclus); D. Gaboriau a montr\'e dans \cite{Gaboriau99} que toutes les valeurs de $[1,\infty]$ \'etaient atteintes.

Un autre invariant de $R$, la \emph{constante isoperim\'etrique} $h(R)$, peut \^etre d\'efini \`a partir des $\Sigma_\Phi$:
\[
h(R)=\inf_{\Phi}\, h(\Sigma_\Phi)
\]
 o\`u $h(\Sigma_\Phi)$ est la constante isop\'erim\'etrique de $\Sigma_\Phi$, d\'efinie de fa\c con appropri\'ee (une d\'efinition pr\'ecise est rappel\'ee \`a la section  \ref{not}). Cet invariant, \`a valeurs dans $[0, \infty]$, a \'et\'e introduit par R. Lyons et les auteurs dans \cite{lpv}. 
 
Dans \cite{lpv} il est montr\'e que
 \[
C(R)\leqslant 1+{1\over 2}h(R),
\] 
si $R$ est de co\^ut fini, et que $h(R)=0$ d\`es que $C(R)=1$. Nous renvoyons \`a cet article pour des d\'etails sur ces deux \'enonc\'es ainsi que  leurs preuves. 
 Le th\'eor\`eme principal de cet article compl\`ete ces deux r\'esultats :

\begin{theorem}\label{th1}
Soit $R$ une  relation d'\'equivalence mesur\'ee ergodique  de type $\IIi$ sur un espace de probabilit\'e $(X,\mu)$ sans atome. Si $C(R)<\infty$, alors 
\[
C(R)=1+{1\over 2}h(R),
\] 
o\`u $C(R)$ est le co\^ut de $R$ et $h(R)$ la constante isop\'erim\'etrique.
\end{theorem}

Le th\'eor\`eme \ref{th1} \'etablit l'\'egalit\'e entre deux invariants de $R$ de nature tr\`es diff\'erente :  
  le co\^ut de $R$ est un invariant \emph{local} des graphes de Cayley de $R$ (par d\'efinition), alors que la constante isop\'erim\'etrique en est un invariant \emph{asymptotique}.
L'\'enonc\'e correspondant en th\'eorie des groupes n'est bien s\^ur pas correct. Si  $\G$ est un groupe d\'enombrable infini de type fini, l'analogue du co\^ut  est le nombre minimal $g(\G)$ de g\'en\'erateurs  de $\G$ (i.e. la valence minimale des graphes de Cayley divis\'ee par 2), et la constante isop\'erim\'etrique de $\G$ est l'infimum $h(\G)$, sur les graphes de Cayley $Y$ de type fini de $\G$, de la constante de isop\'erim\'etrique habituelle de $Y$. En g\'en\'eral, les deux constantes $g(\G)$ et $h(\G)$ ne co\"\i ncident pas modulo renormalisations. Par exemple pour les groupes ab\'eliens libres $\ZI^n$, $n\geqslant 2$, ou plus g\'en\'eralement pour les groupes moyennables non cycliques, on a $h(\G)=0$ et $g(\G)>1$. Une \'egalit\'e analogue \`a celle du th\'eor\`eme \ref{th1} a lieu cependant pour le groupe libre $\G=F_k$  \`a $k$ g\'en\'erateurs; dans ce cas  $g(F_k)=k$ et $h(F_k)=2k-2$. 
En fait, il est facile de voir que pour un groupe de type fini $\G$, l'\'egalit\'e $g(\G)=1+{1\over 2}h(\G)$ n'a lieu que si $\G$ est un groupe libre. En effet, $g(\G)$ est atteint pour un certain syst\`eme g\'en\'erateur (fini) $S$, et  pour ce $S$ on a donc $2(|S|-1)=h(\G)\leqslant h(Y)$, o\`u $h(Y)$ est la constante isop\'erim\'etrique du graphe de Cayley $Y$ de $\G$ relativement \`a $S$. En appliquant cette in\'egalit\'e aux boules $B_n$ de $Y$, on obtient par un calcul imm\'ediat que $|B_n|\geqslant 2|S|(2|S|-1)^{n-1}$ pour tout $n\geq 1$, et donc que $\G$ est un groupe libre. Au contraire, le co\^ut $C(R)$ de $R$ n'est jamais atteint, sauf si $R$ est arborable (cf. \cite{Gaboriau99}).

\subsection*{Remerciements.} Nous remercions George Skandalis qui a sugg\'er\'e que l'in\'egalit\'e de \cite{lpv} \'etait une \'egalit\'e. Le premier auteur est financ\'e par JSPS et l'Universit\'e de Tokyo (IPMU).

\section{Rappels et notations}\label{not}

Faisons quelques brefs rappels pour fixer les notations. Pour plus de d\'etails sur les graphages et le coût des relations d'équivalences, nous renvoyons aux articles \cite{Gaboriau99, Gaboriau02} de D. Gaboriau.

Soit $\Phi$ une famille d\'enombrable d'isomorphismes partiels de $R$, i.e. d'isomorphismes partiels $\varphi$ de $X$ tels que $x\sim \varphi(x)~\mathrm{mod} (R)$. Le \emph{co\^ut} de $\Phi$ est le nombre positif (éventuellement infini)
\[
C(\Phi)=\sum_{\varphi\in \Phi} \mu(\dom \varphi) =\sum_{\varphi\in \Phi} \mu(\Im \varphi)
\] 
o\`u $\dom \varphi$ et $\Im \varphi$ sont respectivement le domaine et l'image de  $\varphi$. On dit que $\Phi$ est un \emph{syst\`eme de g\'en\'erateurs} de $R$ si, presque surement, \'etant donn\'es deux points  $x, y\in X$ \'equivalents modulo $R$, il existe une compos\'ee $\psi=\varphi_n^{\varepsilon_n} \ldots \varphi_1^{\varepsilon_1}$ telle que $\psi(x)=y$, o\`u $\varphi_i\in \Phi$ et $\varepsilon_i\in \{-1,1\}$. Le \emph{co\^ut} de $R$ est l'infimum
\[
C(R)=\inf_\Phi C(\Phi)
\]
pris sur tous les syst\`emes g\'en\'erateurs $\Phi$ de $R$ (voir \cite{Gaboriau99}).
\`A un syst\`eme g\'en\'erateur $\Phi$ de $R$ est naturellement associ\'e un \emph{graphe de Cayley} $\Sigma_\Phi$  de  $R$. Rappelons simplement que $\Sigma_\Phi$ est fibré au-dessus de $X$, et que la fibre $\Sigma_\Phi^x$ en $x\in X$ est un graphe connexe dont les sommets sont les \'el\'ements de $R^x$ et les ar\^etes l'ensemble des couples $\big((x, y), (x, \varphi(y)\big) \in R^x \times R^x$ pour $\varphi\in \Phi$ tel que $y\in \dom \varphi$.
Il est clair que la d\'efinition du co\^ut de $R$ donn\'e ci-dessus co\"\i ncide avec celle donn\'ee en introduction. 

Dans \cite{lpv} nous d\'efinissons une constante isop\'erim\'etrique $h(\Sigma_\Phi)$ de $\Sigma_\Phi$ comme suit:
\[
h(\Sigma_\Phi)=\inf_{A} {{{\nu^{(1)}(\partial_{\Sigma_\Phi} A)}} \over |A|}
\]
o\`u $\nu^{(1)}$ est la mesure de d\'ecompte sur le 1-squelette de $\Sigma_\Phi$, et $A$ un ensemble fini de sommets sym\'etriques de $\Sigma_\Phi$ deux \`a deux disjoints .   Par d\'efinition, on appelle sommet sym\'etrique de $\Sigma_\Phi$ (le graphe d') un automorphisme du groupe plein de $R$ (i.e. un automorphisme de $X$ dont le graphe est inclus dans $R$), et on dit que deux sommets sont disjoints si l'intersection de leurs graphes est n\'egligeable pour la mesure sur $R$ (voir \cite[Definition 4.1]{lpv}). Le bord $\partial_{\Sigma_\Phi} A$ de l'ensemble $A$ est l'ensemble des ar\^etes de $\Sigma_\Phi$ avec exactement un sommet dans $A$ (i.e. dans le graphe d'un automorphisme de $A$). Notons que cette d\'efinition de $h(\Sigma_\Phi)$ ne revient pas \`a \og\ moyenner la constante isop\'erim\'etrique des fibres\fg.

La \emph{constante isop\'erim\'etrique de $R$} est l'infimum
\[
h(R)=\inf_\Phi h(\Sigma_\Phi)
\]
pris  sur tous les syst\`emes g\'en\'erateurs $\Phi$ de $R$.

\section{D\'emonstration du th\'eor\`eme \ref{th1}}\label{s2}

Soit $R$ une relation d'\'equivalence ergodique de type $\IIi$ sur un espace de probabilit\'e $(X,\mu)$ sans atome. Nous montrons ici que
\[
C(R)\geqslant 1+ {1\over 2}h(R),
\] 
l'in\'egalit\'e inverse faisant l'objet de \cite[Theorem 12]{lpv}. On suppose $C(R)<\infty$.
 
Soit $\varphi$ un automorphisme ergodique du groupe plein de $R$. \'Etant donn\'e $\varepsilon>0$, on peut compl\'eter $\varphi$ par une famille $\Phi_\varepsilon'$ d'isomorphismes partiels de $R$ de sorte que 
\[
\Phi_\varepsilon=\Phi_\varepsilon'\cup\{\varphi\}
\]
soit un graphage de $R$, satisfaisant
\[
C(\Phi_\varepsilon )<C(R)+\varepsilon/4
\]
(voir  \cite[Lemme III.5]{Gaboriau99}). 

Soit  $f$ la fonction d\'efinie pour $x\in X$ par
\[
f(x)=\sum_{\psi\in \Phi_\varepsilon'} \chi_{\dom \psi}(x)+\chi_{\Im\psi}(x).
\] 
Elle a pour int\'egrale
\[
\int_X f\d\mu =2C(\Phi_\varepsilon').
\]
Comme la mesure $\mu$ est invariante, il en r\'esulte que
\[
{1\over n+1} \int_X \sum_{i=0}^n f(\varphi^i(x))\d\mu(x) =2C(\Phi_\varepsilon')
\]
pour tout nombre entier $n\geqslant 0$.

Fixons $n> 4/\varepsilon$ et consid\'erons  l'ensemble $A_n$ d\'efini par 
\[
A_n=\{\varphi^i\}_{i=0}^n.
\]
Clairement $A_n$ est un ensemble de points sym\'etriques   deux \`a deux disjoints de $\Sigma_\varepsilon$, au sens de \cite[Definition 8]{lpv}.   De plus pour $x\in X$ le bord de  $A_n^x$ dans $\Sigma_{\varepsilon}^x$ v\'erifie 
\[
|\del_{\Sigma_\varepsilon^x} A^x_n|\leqslant \sum_{i=0}^n f(\varphi^i(x)) + 2,
\]
o\`u $|\cdot|$ est le cardinal. Notons qu'on a \'egalit\'e si $x\not \sim \psi(x)~ \mathrm{mod}~ (R_\varphi)$ pour tout $\psi\in \Phi_\varepsilon'$ et tout $x\in \dom \psi$ --- ce qu'on peut supposer a priori quitte \`a r\'eduire les domaines des \'el\'ements de $\Phi_\varepsilon'$, mais ce n'est pas n\'ecessaire. 

Par suite on a:
\[
{1\over n+1} {\nu^{(1)}(\del_{\Sigma_\varepsilon} A_n) }\leqslant 2C(\Phi_\varepsilon') + {2\over n+1} < 2(C(R)-1)+\varepsilon,
\]
o\`u $\nu^{(1)}$ est la measure naturelle sur le 1-squelette de $\Sigma_\varepsilon$ (associ\'ee \`a la mesure de Haar sur $R$). Comme 
\[
h(R)\leqslant \inf_ A  {{\nu^{(1)}(\del_{\Sigma_\varepsilon} A) }\over |A|}\leqslant {1\over n+1} {\nu^{(1)}(\del_{\Sigma_\varepsilon} A_n) },
\]
on obtient 
\[
h(R)<2(C(R)-1)+\varepsilon,
\]
et lorsque $\varepsilon\to 0$,
\[
C(R)\geqslant 1+ {1\over 2}h(R),
\] 
ce qu'il fallait d\'emontrer.

\section{Commentaires}\label{s3}

\subsection{Formule de compression} Dans \cite{lpv} on montre que le groupe fondamental d'une relation uniform\'ement non moyennable ergodique $R$ est trivial, en \'etablissant l'in\'egalit\'e
\[
h(R)\leqslant \mu(D)h(R_{|D})
\]
o\`u $D$ est une partie bor\'elienne non n\'egligeable de $X$ et $\R_{|D}$ est la restriction de $R$ \`a $D$. En fait cette in\'egalit\'e  est  une \'egalit\'e :

\begin{corollary} Soit $R$ comme dans le th\'eor\`eme \ref{th1} et $D$ une partie bor\'elienne non n\'egligeable de $X$. On a:
\[
h(R)= \mu(D)h(R_{|D}).
\]
\end{corollary}

Ceci r\'esulte imm\'ediatement du th\'eor\`eme \ref{th1} et la proposition II.6 de \cite{Gaboriau99} concernant  le co\^ut.

\subsection{Une autre d\'efinition du bord} Conservant les notations introduites ci-dessus, on peut aussi d\'efinir le bord $\del_\Sigma A$ comme l'ensemble des \emph{points de $R$} dans le bord de $A$ relativement \`a $\Sigma$ (plut\^ot que les \emph{ar\^etes de $\Sigma$} dans le bord de $A$ relativement \`a $\Sigma$). On obtient alors une constante $h'(R)$, a priori plus petite que $h(R)$, et les r\'esultats de \cite{lpv} \'etablissent l'in\'egalit\'e (voir la note p. 598) 
\[
C(R)\leqslant 1+{h'(R)\over 2}.
\]
On a donc en fait, d'apr\`es le th\'eor\`eme \ref{th1}:
\[
C(R)= 1+{h(R)\over 2}= 1+{h'(R)\over 2}.
\]
Comme not\'e dans \cite[p. 598]{lpv}, la version $h'$ de la constante isop\'erim\'etrique est \'evidemment major\'ee par le taux de croissance uniforme $\omega(\G)$ de $\G$ moins 1, 
et on a donc   
\[
\omega(\G)\geqslant 2\beta_1(\G)+1
\]
(voir \cite[Cor. 3.2]{lpv}), o\`u $\beta_1(\G)$ est le premier nombre de Betti $\ell^2$ de $\G$.    

\subsection{Taux de cycles} Rappelons qu'\`a une relation d'\'equivalence $R$ ergodique de type $\IIi$ est associ\'ee une suite de nombres positifs $\beta_0(R),\beta_1(R),\beta_2(R),\ldots$, appel\'es nombres de Betti $\ell^2$ de $R$ (voir \cite{Gaboriau02}).  On a toujours
\[
 C(R)\geqslant\beta_1(R)+1
 \]
et on se sait pas si l'in\'egalit\'e peut \^etre stricte (voir la section 3.6 de \cite{Gaboriau02}).
 Le \emph{taux de cycle} de $R$ (voir \cite{th})  est d\'efini comme suit: 
 \[
 \tau(R)=\inf_{\Phi} \dim_{LR} \overline{Z_1(\Sigma_\Phi)},
 \]
o\`u l'infimum est pris sur les syst\`emes g\'en\'erateurs $\Phi$ de $R$,  $\Sigma_\Phi$ est le graphe de Cayley de $\Phi$, $Z_1(\Sigma_\Phi)$ est l'espace des cycles de $\Sigma_\Phi$. La dimension $\dim_{LR}$ est la dimension de l'adh\'erence Hilbertienne $\overline{Z_1(\Sigma_\Phi)}$  prise au sens de Murray et von Neumann (comme pour les nombres de Betti $\ell^2$), où $LR$ désigne l'alg\`ebre de von Neumann de $R$. On a alors
 \[
 C(R)=\tau(R)+\beta_1(R)+1
 \]
pour toute relation ergodique  $\IIi$ de co\^ut fini (voir \cite{th,sc}), et
 il n'y a pas d'exemples connus pour lesquels $\tau(R)>0$. Le taux de cycles a \'et\'e introduit dans \cite{th} conjointement avec une technique \emph{d'effa\c cage de cycles}, qui permet de faire diminuer la dimension $\dim_{LR} \overline{Z_1(\Sigma_\Phi)}$.  En fait (cf. \cite{th}) il est possible \emph{d'effacer r\'ecursivement  tous les cycles} d'un graphage donn\'e---mais ceci ne suffit pas a priori \`a entra\^\i ner $\tau(R)=0$. Le th\'eor\`eme \ref{th1} pourrait fournir une nouvelle approche \`a ces questions.

\end{document}